\documentclass{birkjour_t2}
 \usepackage{graphicx, subfigure}
 
\usepackage{amssymb, amsmath, latexsym}
\usepackage{amsfonts}
\allowdisplaybreaks

\catcode`\@=11
\@addtoreset{equation}{section}
\def\theequation{\thesection.\arabic{equation}}
\def\theequation{\thesection.\arabic{equation}}
\newcommand \RR   {\mathbb{R}}
\newcommand \del  {\partial}  
  
\newcommand \Ga   {\Gamma}

\newcommand \be   {\begin{equation}}
\newcommand \ee   {\end{equation}}

 \newtheorem{theorem}{Theorem}[section]
 
 \newtheorem{lemma}[theorem]{Lemma}
 
 \theoremstyle{definition}
 \newtheorem{definition}[theorem]{Definition}
 \theoremstyle{remark}
 \newtheorem{remark}[theorem]{Remark}
 
 \numberwithin{equation}{section}

\newcommand{\bp}{\begin{proof}}
\newcommand{\ep}{\end{proof}}

\begin{document}

\title[Isometric Immersion with Slowly Decaying Gauss Curvature]{Isometric Immersions via Compensated Compactness for\\ Slowly Decaying Negative Gauss Curvature and Rough Data}

\author{Cleopatra Christoforou} 
\address{C.C: Department of Mathematics and Statistics\\ University of Cyprus\\ 1678 Nicosia\\ Cyprus}
\urladdr{http://www2.ucy.ac.cy/$\sim$kleopatr}
\email{Christoforou.Cleopatra@ucy.ac.cy}

\author{Marshall Slemrod}
\address{M.S: Department of Mathematics\\ University of Wisconsin\\ Madison, USA}
\email{slemrod@math.wisc.edu}


\begin{abstract}
In this paper the method of compensated compactness is applied to the problem of isometric immersion of a two dimensional Riemannian manifold with negative Gauss curvature into three dimensional Euclidean space. Previous applications of the method to this problem have required decay of order $t^{-4}$ in the Gauss curvature. Here we show that the decay of Hong~\cite{H} $t^{-2-\delta/2}$ where $\delta\in(0,4)$ suffices.
\end{abstract}

\maketitle

\section{Introduction}\label{S1}
In two recent papers Chen, Slemrod and Wang~\cite{CSW} and Cao, Huang and Wang~\cite{CHW} have used the method of compensated compactness to establish global isometric immersions into $\RR^3$ for two dimensional Riemannian manifolds for rough data. In the examples considered the Gauss curvature was negative and decayed at least as $t^{-4}$ where initial data was given at $t=0$. Needless to say that leaves open the question as whether the compensated compactness method will work for a slower rate of decay. Certainly based on the paper of Hong~\cite{H} which has also been exposited in the book of Han and Hong~\cite{HH} we expect the result to be true for decay of order $t^{-2-\delta/2}$ where $\delta$ is between $0$ and $4$. The proof of Hong is a careful study of the hyperbolic system of two balance laws (the two Codazzi equations) and one closure relation (the Gauss equation) and requires two separate steps. The first step is to establish existence of smooth solutions to the balance laws for small, smooth data prescribed at a large enough time $t=T_1$. The reason for this part is that is only after large time that the decay of the Gauss curvature may be exploited to obtain the relevant $C^1$ a priori estimates. The second part of the proof is rather standard and simply  asks for the initial data at $t=0$ to be sufficiently small and smooth to enable us to get a solution up to $t=T_1$. Here no reference is given in Han and Hong but a standard existence, uniqueness theorem for quasi-linear hyperbolic systems will suffice. Such a theorem may be found in Janenko and Rozdestvenskii~\cite[Chapter 1, Sect. 8]{JR} where the growth in $C^1$ of solutions is governed by a coupled pair of ordinary differential equations, one which is of Ricatti type. Hence just as in the classical theory of ordinary differential equations small data allows for a longer time of existence. In this paper as compliment to Hong's result we reconsider the first part of Hong's program and show that in fact that rough $L^\infty$ data suffices at the initial time $t=T_1$ and that the method of compensated compactness will yield existence of weak solutions to the Gauss-Codazzi system for $t>T_1$. Of course this does not provide a new result since for the second part of Hong's proof one would still need the classical smooth solution existence, uniqueness theorem to reach time $t=T_1$ from $t=0$. Nevertheless we believe that the new application of the compensated compactness method is appealing and is of independent interest.

Before continuing a short historical view is in order. First we note the relevant Gauss-Codazzi system can be written as a linearly degenerate system or what is termed ``weakly non-linear quasi-linear" system in the monograph of Janenko and Rozdestvenskii~\cite{JR} and is discussed in Chap. $4$, Sec. $4$ of that book. They note that such systems possess the property that uniform boundedness of solutions on $0\le t\le T$ and strict hyperbolicity imply uniform boundedness of first derivatives on $0\le t\le T$ if these derivatives are initially bounded. Thus it appears that the crucial estimates will be a uniform bound on the dependent variables and in addition proof that the strict hyperbolicity is not lost. It was this path that was followed by Hong~\cite{H} and yields the following result:

\begin{theorem}[Han and Hong~{\cite[Theorem 10.2.2]{HH}}]\label{HH}
For a complete simply
connected two dimensional Riemannian manifold $(\mathcal{M},g)$ with negative Gauss curvature $K$
with metric $g=dt ^{2}+B^{2}(x,t) dx^{2}$. Assume for some constant $\delta>0$
\begin{enumerate}
\item[(i)] $t ^{2+\delta }\left\vert K\right\vert $ is decreasing in $|t|$, $|t|>T$;
\item[(ii)] $\partial _{x }^{i}\ln \left\vert K\right\vert$, for $i=1,2$ and $t \partial
_{t }\partial _{x }\ln \left\vert K\right\vert $ are bounded;
\item[(iii)] $K$ is periodic in $x$ with period $2\pi$.
\end{enumerate}
Then $(\mathcal{M},g)$ admits a smooth isometric immersion in $
\mathbb{R}
^{3}.$
\end{theorem}

But as noted above a search for ``corrugated immersions" would ask that the data be ``rough" and not in $C^1$  and that is the issue pursued here. In particular, we study data in $L^\infty$. It should be mentioned that for discontinuous data of bounded variation, isometric immersions have been established using a different method in~\cite{Cisom} but again with decay rate at least as $t^{-4}$. An exposition of the current state of the theory of systems of balance laws can be found in the book~\cite{Dafermos3}.

The paper contains six sections after this introduction. Section~\ref{S2} provides a review of the isometric embedding problem and exposits the theorem of S. Mardare on non-smooth embeddings. In Section~\ref{SV} we give a viscous approximation scheme for resolving the relevant balance laws: the equations of Gauss and Codazzi. In Section~\ref{Sbound} we derive apriori $L^\infty$ estimates for the viscous approximations
and in Section~\ref{SH} we show the viscous system possesses a crucial $H^{-1}_{loc}$ estimate which is needed to apply the method of compensated compactness.  In Section~\ref{SCC} we recall the compensated compactness framework~\cite{CSW, CSW2} and show that passage to the inviscid limit may be accomplished hence yielding the desired non--smooth immersion. Finally Section~\ref{S7} pursues the issue as to whether the decay rate of the Gauss curvature can be reduced, say as given by the choice
$$ K=-\frac{1}{(3+t)^2(\ln(3+t))^p},\qquad t>0,\,p \text{\,\,sufficiently large}.$$ Here we show that for this choice one of the key a priori estimates -- the preservation of strict hyperbolicity -- is retained. Hence any lack of non-smooth embedding  must be due to lack of $L^\infty$ bounds on two of the three components on the second fundamental form (the third component is a priori bounded).

\section{Preliminaries}\label{S2}
Let $\Omega\subset\RR^2$ be an open set. Consider a map ${\bf y}:\Omega\to\RR^3$ having the tangent plane of the surface ${\bf y}(\Omega)\subset\RR^3$ at ${\bf y}(x_1,x_2)$ spanned  by the vectors $\{\del_1 {\bf y},\del_2 {\bf y}\}$. Then, the unit normal vector ${\bf n}$ to the surface ${\bf y}(\Omega)$ is given by
\be\label{S1:normal}
	{\bf n}=\frac{\del_1 {\bf y}\times\del_2 {\bf y}}{|\del_1 {\bf y}\times\del_2 {\bf y}|}
\ee
and the corresponding metric is
$$
ds^2=d{\bf y}\cdot d{\bf y}\;,
$$
or equivalently,
\be
ds^2= (\del_1 {\bf y}\cdot\del_1 {\bf y})(dx_1)^2+2(\del_1 {\bf y}\cdot\del_2 {\bf y})dx_1\,dx_2+(\del_2 {\bf y}\cdot\del_1 {\bf y})(dx_2)^2\;.
\ee
The isometric immersion problem is an \emph{inverse} problem: Given $(g_{ij})$ $i,j=1,2$ functions in $\Omega$, with $g_{12}=g_{21}$, find a map ${\bf y}:\Omega\to\RR^3$ so that
\be
 	d{\bf y}\cdot d{\bf y}=g_{11} (dx_1)^2+2g_{12} dx_1dx_2+g_{22}(dx_2)^2\;,
\ee
or equivalently,
\be
\del_1{\bf y}\cdot\del _1{\bf y}=g_{11},\qquad \del_1{\bf y}\cdot\del _2{\bf y}=g_{12},\qquad \del_2{\bf y}\cdot\del _2{\bf y}=g_{22}
\ee
with a linearly independent set $\{\del_1 {\bf y},\del_2 {\bf y}\}$ in $\RR^3$. Hence, the isometric immersion problem is fully nonlinear in the three unknowns being the three components of the map ${\bf y}$.

We recall that a two dimensional manifold $(\mathcal{M},g)$ parametrized by $\Omega$ with associated metric $g=(g_{ij})$ admits two fundamental forms: the first fundamental form $I$ for $\mathcal{M}$ on $\Omega$ is
\be
	I\doteq g_{11} (dx_1)^2+2g_{12} dx_1dx_2+g_{22}(dx_2)^2
\ee
and the second fundamental form $II$ is
\be
 II\doteq -d{\bf n}\cdot d{\bf y}=h_{11}(dx_1)^2+2h_{12} dx_1\,dx_2+h_{22}(dx_2)^2
\ee
with ${\bf n}$ being the unit normal vector to $\mathcal{M}$. The coefficients $(h_{ij})$ represent the orthogonality of ${\bf n}$ to the tangent plane and are associated with the second derivatives of ${\bf y}$ and
since ${\bf n}\cdot d{\bf y}=0$, it follows
$$
 II=({\bf n}\cdot \del_1^2{\bf y})(dx_1)^2+ ({\bf n}\cdot\del_1\del_2{\bf y})dx_1\,dx_2+ ({\bf n}\cdot\del_2^2{\bf y})(dx_2)^2\;.
$$

By equating the cross-partial derivatives of ${\bf y}$, the isometric immersion problem as stated above reduces to the Gauss--Codazzi system
\be\label{S2:GC}
\begin{array}{l}
\del_1 M-\del_2 L=\Ga_{22}^{(2)}L-2\Ga_{12}^{(2)}M+\Ga_{11}^{(2)} N\\
\del_1 N-\del_2 M=-\Ga_{22}^{(1)}L+2\Ga_{12}^{(1)}M-\Ga_{11}^{(1)} N
\end{array}
\ee
with the condition
\be\label{S2:GCrel}
 LN-M^2=K\,,
\ee
where
\be\label{S2:hrelation}
 L=\frac{h_{11}}{\sqrt{|g|}},\qquad  M=\frac{h_{12}}{\sqrt{|g|}}, \qquad  N=\frac{h_{22}}{\sqrt{|g|}}
\ee
and
$|g|\doteq det(g_{ij})=g_{11} g_{22}-g^2_{12}$. The Gauss curvature $K=K(x_1,x_2)$ is given by
\be
	K(x_1,x_2)=\frac{R_{1212}}{|g|}\;,
\ee
where $R_{ijkl}$ is the curvature tensor 
\be
	R_{ijkl}=g_{lm}\left(\del_k\Ga_{ij}^{(m)}-\del_j\Ga_{ik}^{(m)}+\Ga_{ij}^{(n)}\Ga_{nk}^{(m)}-\Ga_{ik}^{(n)}\Ga_{nj}^{(m)}\right)\;,
\ee
and $\Ga_{ij}^{(k)}$ is the Christoffel symbol
\be
	\Ga_{ij}^{(k)}\doteq \frac{1}{2} g^{kl}\left( \del_j g_{il}+\del_i g_{jl}-\del_lg_{ij}\right)\;.
\ee
Here, the indices $i,\,j,\,k,\,l=1,2$, $(\del_1,\del_2)=(\del_{x_1},\del_{x_2})$ and the summation convention is used. Also, $(g^{kl})$ is the inverse of $(g_{ij})$. 

The fundamental theorem of surface theory states that given forms $I$ and $II$ with $(g_{ij})$ being positive definite and smooth coefficients, $(g_{ij})$ and $(h_{ij})$ that satisfy the Gauss-Codazzi system~\eqref{S2:GC}--\eqref{S2:hrelation}, then there exists a surface embedded into $\RR^3$ with first and second fundamental forms $I$ and $II$. This result has been extended by S.~Mardare~\cite{M1} when  $(h_{ij})\in L^\infty_{loc}(\Omega)$ for given $(g_{ij})\in W^{1,\infty}_{loc}(\Omega)$ and then, the surface immersed is $C^{1,1}(\Omega)$ locally. Thus, the isometric immersion problem reduces to solving the Gauss-Codazzi system~\eqref{S2:GC}--\eqref{S2:hrelation} for $(h_{ij})\in L^{\infty}_{loc}(\Omega)$ with a given positive definite metric $(g_{ij})\in W^{1,\infty}_{loc}(\Omega)$ and then, immediately, we recover the immersion surface ${\bf y}({\Omega})$, which is  $C^{1,1}$ locally. We refer the reader to books~\cite{DoC,HH} for an exposition of the surface theory and to Mardare~\cite{M1,M2} for the extension of the aforementioned result to $(h_{ij})\in L^\infty_{loc}$.

However for completeness of our presentation we state Mardare's result in full:
\begin{theorem}[S.~Mardare~\cite{M1}]\label{Mardare}
Assume that $\Omega$ is a connected and simply-connected open subset of $\RR^2$ and that the matrix fields $(g_{ij})\in W^{1,\infty}_{loc}(\Omega)$ being symmetric positive definite and $(h_{ij})\in L^\infty_{loc}(\Omega)$ symmetric satisfy the Gauss and Codazzi-Mainardi equations in $\mathcal{D}'(\Omega)$. Then there exists a mapping ${\bf y}\in W^{2,\infty}_{loc}(\Omega,\RR^3)$ such that
\begin{align}
g_{ij}&=\partial_i {\bf y}\cdot\partial_j {\bf y},\nonumber\\
h_{ij}&=\partial_{ij} {\bf y}\cdot\frac{\partial_1{\bf y}\times \partial_2 {\bf y}}{|\partial_1{\bf y} \times\partial_2 {\bf y}|}\nonumber
\end{align}
a.e. in $\Omega$. Moreover, the mapping ${\bf y}$ is unique in $W^{2,\infty}_{loc}(\Omega,\RR^3)$ up to proper isometries in $\RR^3$.
\end{theorem}

By virtue of the embedding of $W^{2,\infty}$ into $C^{1,1}$ the immersion ${\bf y}$  is locally in $C^{1,1}$, cf. Evans~\cite[Chapter 5]{Ev}. Now let us recall the following definition.

\begin{definition} $(\mathcal{M},g)$ is a geodesically complete Riemannian manifold
if and only if every geodesic can be extended indefinitely.
\end{definition}

It is perhaps useful to note that the concept of geodesically complete Riemannian manifold is equivalent to the Riemannian manifold defining a complete metric space. This is a consequence of Hopf--Rinow theorem, cf.~\cite[Chapter 7]{DoC}. In fact, do~Carmo notes ``intuitively, this means that the manifold does not have any holes or boundaries".

Under the assumption that our two dimensional manifold is geodesically
complete and simply connected  we can simplify the structure of our metric.
The exact result is as follows and is essentially due to Hadamard but we use
the presentation given in Han and Hong~\cite{HH}.

\begin{lemma}[Han and Hong~{\cite[Lemma 10.2.1]{HH}}] Let $(\mathcal{M},g)$ be a geodesically complete simply connected smooth two
dimensional Riemannian manifold with non-positive Gauss curvature. Then
there exists a global geodesic coordinate system $(x,t)$ in $\mathcal{M}$ with metric 
\be\label{S2 g} g=dt^{2}+B^{2}(x,t)dx^{2}\ee
where $B$ is a smooth function satisfying $
B(x,0)=1$ and $\partial _{t}B(x,0)=0$ for $x\in \RR$.
\end{lemma}
A direct substitution of~\eqref{S2 g} in~\eqref{S2:GC} then yields that $L,M,N$ satisfy the Gauss-Codazzi
system in the form
\be\label{S2:GC new}
\begin{array}{l}
\partial _{t}L-\partial _{x}M=L\partial _{t}\ln B-M\partial _{x}\ln B+NB\partial _{t}B,\\
\partial _{t}M-\partial _{x}N=-M\partial _{t}\ln B,
\end{array}
\ee
with
\be LN-M^{2}=KB^{2},\ee
where $\partial _{tt}B=-KB$ defines the Gauss curvature $K$ in terms of the
metric.

\section{The Viscous Approximation}\label{SV}
As we will be dealing with non--smooth data a natural approach is embed our initial value problem in viscous approximating system with viscosity $\mu>0$ and attempt to recover our solution as limit for $\mu\to0+$. More precisely, in this section, we first study the viscous approximation of the scaled variables of $(L,M,N)$ to system~\eqref{S2:GC new} and then establish properties for a family of metrics of the form~\eqref{S2 g} that correspond to the class of negative curvature with decay rate of the order of $t^{-2-\frac{\delta}{2}}$ with $\delta\in(0,4)$.

First it is easy to check that the scaled variables
\be l=\frac{L}{B^2\sqrt{|K|}},\quad m=\frac{M}{B\sqrt{|K|}},\quad n=\frac{N}{\sqrt{|K|}}\ee
satisfy the system
\be\label{SV:lm}
\begin{array}{l}
\partial _{t}l-\frac{1}{B}\partial _{x}m+(l-n)\partial _{t}\ln B+\frac{l}{2}\partial _{t}\ln |K|  -\frac{m}{2B}\partial _{x}\ln |K|=0,\\\\
\partial _{t}m-\frac{1}{B}\partial _{x}n+2m\partial _{t}\ln B+\frac{m}{2} \partial _{t}\ln |K| -\frac{n}{2}\partial _{x}\ln |K|   =0,
\end{array}
\ee
with
\be\label{SV:rel} ln-m^{2}=-1. \ee
The eigenvalues associated with system~\eqref{SV:lm} are
\be\lambda_1=\frac{m-1}{l},\qquad\lambda_2=\frac{m+1}{l}\ee 
and we see that each characteristic field is linear degenerate. System~\eqref{SV:lm} is \emph{strictly hyperbolic} if $\lambda_1<\lambda_2$, or equivalently if $l$ is finite.

Consider the viscous approximations $(l^\mu, m^\mu, n^\mu)$ that satisfy system
\be\label{SV:lm viscous}
\begin{array}{l}
\partial _{t}l^\mu-\dfrac{1}{B}\partial _{x}m^\mu+(l^\mu-n^\mu)\partial _{t}\ln B+\dfrac{l^\mu}{2}\partial _{t}\ln |K|  -\dfrac{m^\mu}{2B}\partial _{x}\ln |K|=\mu\partial _{xx}l^\mu,\\\\
\partial _{t}m^\mu-\dfrac{1}{B}\partial _{x}n^\mu+2m^\mu\partial _{t}\ln B+\dfrac{m^\mu}{2} \partial _{t}\ln |K| -\dfrac{n^\mu}{2}\partial _{x}\ln |K|   =\mu\partial _{xx}m^\mu,
\end{array}
\ee
with
\be\label{SV:rel viscous} l^\mu n^\mu-(m^\mu)^{2}=-1. \ee
Here $\mu>0$ is a constant ``viscosity". For convenience we drop the index $\mu$ from the viscous approximate solutions $(l^\mu, m^\mu, n^\mu)$ from here and on and we use it again in Section~\ref{SCC} when studying the limit $\mu\to 0+$.

Set the Riemann invariants 
\be u=-\frac{m}{l}+\frac{1}{l},\qquad v=-\frac{m}{l}-\frac{1}{l}\ee

\be\label{SV: lmn} l=\frac{2}{u-v},\quad m=-\left(\frac{u+v}{u-v}\right),\quad n=\frac{2uv}{u-v}\ee

Multiply system~\eqref{SV:lm viscous} from the left by $(\partial_l u,\partial_m u)^T$ and $(\partial_l v,\partial_m v)^T$ to obtain the viscous equations of $(u,v)$:

\begin{align}
\partial _{t}u+\frac{v}{B}\partial _{x}u+&v(1+u^2)\partial _{t}\ln B-\frac{(u-v)}{4}\left(\partial _{t}\ln |K|  +\frac{u}{B}\partial _{x}\ln |K|\right)=\nonumber\\
\label{SV:uv viscous1}
&=\frac{\mu}{(v-u)}\left\{ ((\partial_x u)^2-(\partial_x v)^2)-\frac{2u}{v-u} (\partial_x u-\partial_x v)^2-u\partial_{xx} u \right\},\\
\partial _{t}v+\frac{u}{B}\partial _{x}v+&u(1+v^2)\partial _{t}\ln B-\frac{(v-u)}{4}\left(\partial _{t}\ln |K|+\frac{v}{B}\partial _{x}\ln |K|\right)   =\nonumber\\
\label{SV:uv viscous2}
&=\frac{\mu}{(v-u)}\left\{ ((\partial_x u)^2-(\partial_x v)^2)-\frac{2v}{v-u} (\partial_x u-\partial_x v)^2+v\partial_{xx} v \right\}.
\end{align}

It is easy to check that strict hyperbolicity in the $(u,v)$ variables is equivalent to $u\neq v$. Now the system is \emph{uniformly strictly hyperbolic} if $v-u$ is uniformly bounded away from zero.
 
In the following sections, we prove uniform $L^\infty$ bounds independent of $\mu$ to $(u,v)$ and therefore to $(l^\mu,m^\mu, n^\mu)$ via~\eqref{SV: lmn}.

\subsection{The Metric for a Special Case}
Let $(\mathcal{M},g)$ be a geodesically complete simply connected smooth two
dimensional Riemannian manifold with non-positive Gauss curvature $K$ and a metric of the form~\eqref{S2 g}.

To keep our ideas clear and the presentation relatively simple we here consider only the special case given below. Of course the method could be generalized beyond this case at the cost of greater complications and technicalities.

Set $h\doteq B>0$ and $k^*\doteq|K|$ and we assume that $k^*$, $h$ are taken to be independent of $x$. Then $h$ and $k^*$ satisfy 
\be\label{SV: h}
\partial_{tt} h=k^* h,\quad h(0)=1,\quad \partial_{t} h(0)=0.
\ee
Moreover, let $\phi=\phi(t)$ be the solution of 
\be\label{SV phi}
\begin{array}{l}
\partial_{t}\phi=\phi(1+\phi^2)\partial_{t}\ln h+\dfrac{\phi}{2}\partial_{t}\ln k^*,\quad t>T\\\\
\phi(T)=\psi_0   
\end{array}
\ee
with $\psi_0$ a constant. As noted in~\cite[10.2.36]{HH}, the explicit solution $\phi$ of the above problem is given by the expression
\be \phi(t)=\frac{bh\sqrt{k^*}}{(1-2b^2\int_T^t h\,\partial_t h \, k^* ds)^{1/2}}\ee
with 
\be b=\frac{\psi_0}{h(T)\sqrt{k^*(T)}}.\ee
Notice that for $\phi$ to be defined for all $t>T$, we must have $h\partial_t h \, k^*\in L^1(T,\infty)$. By the choice
\be\label{SV: k} k^*=\frac{C}{(1+|t|)^{2+\frac{\delta}{2}}},\quad C>0\ee
as taken in Hong~\cite{H} and $\psi_0$ small enough, then we have the formula
\be\label{SV decr} (1+\phi^2)\partial_t\ln h+\frac{1}{2}\partial_t\ln k^*=-\left(\frac{\delta-4\phi^2}{4t}\right)+O(\frac{1}{t^{1+\delta/2}})\le 0\ee
for $t>T$. Thus, for $\psi_0$ small enough, equation~\eqref{SV phi} implies that $\partial_t\phi<0$, $\phi>0$ for $t>T$. In fact, formula~\eqref{SV decr} is the key to the rest of our analysis.

In what follows, we establish estimates on $h$ and $\partial_t h$ for the chosen decay rate of the curvature~\eqref{SV: k} that are used in the following sections.

\begin{lemma}\label{SV: lemma k}
If $k^*(s)$ and $sk^*(s)$ are in $L^1(0,\infty)$, then
\be\label{SV:la}\int_0^t k^*(s)ds\le \partial_t h\le C_1\ee
and
\be\label{SV:lb} 1+\int_0^t\int_0^s k^*(\tau)dt\tau ds\le h(t)\le 1+C_1 t\ee
\end{lemma}
\bp The proof is given in Han and Hong~\cite[Lemma 10.2.3]{HH}.
\ep

Next, an important estimate for $\partial_t\ln h$ is given in the following lemma.  

\begin{lemma}\label{SV:l2}
Let $|k^* t^{2+\delta}|$ be decreasing in $|t|$ for $|t|>T$. Then
\be \partial_t \ln h=\frac{1}{t}+O(\frac{1}{|t|^{1+\delta}})\ee
for sufficiently large $|t|$ and $\partial_t\ln h$ is bounded.
\end{lemma}
\bp
Again the proof is given in Han and Hong~\cite[Lemma 10.2.3]{HH}.
\ep

\section{Invariant Regions-- $L^\infty$ bounds}\label{Sbound}
In this section we establish $L^\infty$ bounds on the solutions to the viscous system~\eqref{SV:lm viscous}--\eqref{SV:rel viscous}. By~\eqref{SV decr} and the choice of curvature~\eqref{SV: k}, we can get information on the sign of $\partial_t\ln h+\frac{1}{2}\partial_t\ln k^*$. In fact this follows from Lemma~\ref{SV:l2}. Simply write
\be \partial_t\ln h+\frac{1}{2}\partial_t\ln k^*=-\frac{\delta}{4t}+O(\frac{1}{|t|^{1+\delta/2}})<0\ee
for large $t$. Here however we shall use a related equality which we call the sign--switch property. Specifically compute
\begin{align}
\partial_t\ln h+\frac{1}{4}\partial_t\ln k^*&=\frac{1}{t}+O(\frac{1}{|t|^{1+\delta/2}})-\frac{1}{4}(2+\frac{\delta}{2})\frac{1}{t}\nonumber\\
&=\frac{1}{2t}(1-\frac{\delta}{4})+O(\frac{1}{|t|^{1+\delta/2}})
\end{align}

Hence, when $0<\delta<4$, we have the sign-switch
\be\label{Sbound: sign}\partial_t\ln h+\frac{1}{4}\partial_t\ln k^*>0\quad \text{for } t\text{ large}\ee

We now establish the following lemma under the choice $k^*$ given by~\eqref{SV: k} and $0<\delta<4$.

\begin{lemma}\label{Sbound: lemma} Assume that $h$ and $k^*$ are independent of $x$ satisfying~\eqref{SV: h} and $k^*$ is given by~\eqref{SV: k}. Then, if initially at time $t=T_1$, with $T_1$ sufficiently large, the solutions $(u,v)$ of system~\eqref{SV:uv viscous1}--\eqref{SV:uv viscous2} satisfy the bounds
\be\label{Sbound: lbd} -\psi_0\le u\le -e^{-t}\psi_0,\quad e^{-t}\psi_0\le v\le \psi_0\ee
then these bounds on $(u,v)$ persist for all $t>T_1$.
\end{lemma}
\bp
Without loss of generality let us first establish the bound on $u$ from above. If at some first passage time $t^*>T_1$ we have for some $x^*$ that $u(x^*,t^*)=-e^{-t^*} \psi_0$, then at $(x^*,t^*)$ we have $\partial_x u=0$, $\partial_{xx} u\le0$, $\partial_t(u+e^{-t}\psi_0)>0$ and $e^{-t}\psi_0\le v\le \psi_0$. From~\eqref{SV:uv viscous1} we have at $(x^*,t^*)$
\be \partial_t u\le -v\left[(1+e^{-2t^*}\psi_0^2)\partial_t\ln h+\frac{1}{4}\partial_t\ln k^*\right]-\frac{1}{4} e^{-t^*} \psi_0\partial_t\ln k^*\ee
From~\eqref{Sbound: lbd} and~\eqref{Sbound: sign}, we have
\be \partial_t u\le -e^{-t^*} \psi_0 (e^{-2t^*}\psi_0^2 \partial_t\ln h+\frac{1}{4} \partial_t\ln k^*)\ee
for $T_1$ sufficiently large. This is equivalent to 
\be \partial_t (u+e^{-t}\psi_0)\le -e^{-t^*} \psi_0 \left(1+e^{-2t^*}\psi_0^2 \partial_t\ln h-(\frac{1}{2}+\frac{\delta}{8})\frac{1}{1+t^*}\right)\ee
at $(x^*,t^*)$. Hence, for $T_1$ sufficiently large $\partial_t (u+e^{-t}\psi_0)<0$, which is a contradiction of our assumption $\partial_t(u+e^{-t}\psi_0)>0$. A similar argument yields the bound on $v$ from below.

Next we establish the bound on $u$ from below. Again let $t^*$ denote the first passage time. Hence we have some point $(x^*,t^*)$ for which $\partial_x u=0$, $\partial_{xx} u\ge0$, $\partial_t u>0$ with $e^{-t}\psi_0\le v\le \psi_0$. From~\eqref{SV:uv viscous1} we have at $(x^*,t^*)$
\begin{align}
 \partial_t u=&-v(1+\psi_0^2) \partial_t\ln h+\frac{1}{4}(-\psi_0-v)\partial_t\ln k^*\nonumber\\
&+\frac{\mu}{v+\psi_0}\left\{(\partial_x v)^2(-1+\frac{2\psi_0}{v+\psi_0})+\psi_0\partial_{xx} u  \right\}
\end{align}
Adding to this expression the identity~\eqref{SV phi} yields
\begin{align}  \partial_t (u+\phi)&> \phi(1+\phi^2) \partial_t\ln h+\frac{\phi}{2}\partial_t\ln k^*\nonumber\\
&-v(1+\psi_0^2) \partial_t\ln h+\frac{1}{4}(-\psi_0-v)\partial_t\ln k^*
\end{align}
at $(x^*,t^*)$. Now choose the data for $\phi$ to be $\phi(t^*)=\psi_0$. Hence at $(x^*,t^*)$
\be \partial_t (u+\phi)\ge(\psi_0-v)\left[(1+\psi_0^2)\partial_t\ln h+\frac{1}{4}\partial_t\ln k^*\right]\ee
By~\eqref{Sbound: lbd} and the fact that $v\le \psi_0$ we have that $\partial_t (u+\phi)\ge 0$ at $(x^*,t^*)$. Since $u(x^*,t^*)+\phi(t^*)=-\psi_0+\psi_0=0$ we have
$u(x^*,t)\ge-\phi(t)$ for $t>t^*$, $|t-t^*|$ small. Recalling the analysis in~\eqref{SV decr} for the choice~\eqref{SV: k}, $\phi(t)$ is decreasing. Thus we have $u(x^*,t)\ge-\phi(t^*)=-\psi_0$ which is a contradiction. A similar argument holds for $v$. This completes the proof of the lemma.
\ep

\begin{remark}
Of course the proof is motivated by the one given in Han--Hong~\cite{HH}. The advantage of the one given above is its relative simplicity and the precise estimates for $(u,v)$ from above, below (respectively).
\end{remark}
\section{$H^{-1}$ compactness}\label{SH}
In this section, we prove the $H^{-1}$ compactness of the sequence
\be\label{SH: seq} \partial_t l-\partial_x(\frac{m}{h}),\qquad  \partial_t m-\partial_x(\frac{n}{h})\ee
constructed from the viscous approximation~\eqref{SV:lm viscous}--\eqref{SV:rel viscous} with data satisfying~\eqref{Sbound: lbd}. First from the formulas~\eqref{SV: lmn} and the bounds established in Lemma~\ref{Sbound: lemma} we immediately have

\begin{lemma} Under the assumptions of Lemma~\ref{Sbound: lemma}, for $u$, $v$ satisfying~\eqref{Sbound: lbd} at $t=T_1$, $T_1$ sufficiently large, we have $(l,m,n)\in L^\infty([T_1,T_2]\times\RR^3)$ where $T_2$ is any value $T_2>T_1$. 
\end{lemma}

For convenience, let us set $\Omega=[T_1,T_2]\times\RR^3$ so that $(l,m,n)\in L^\infty(\Omega)$. Inspection of~\eqref{SV:lm viscous} shows that $\partial_t l-\partial_x(\frac{m}{h})$,  $ \partial_t m-\partial_x(\frac{n}{h})$ will lie in a compact subset of $H^{-1}_{loc}(\Omega)$ if the viscous terms $\mu\partial_{xx} l$, $\mu\partial_{xx} m$ will lie in a compact subset of $H^{-1}_{loc}(\Omega)$. To show this we follow a standard argument say as given in the paper of Cao, Huang and Wang~\cite{CHW}.

Define the entropy, entropy flux pair 
\be \eta=-\frac{m^2+1}{l},\qquad q=\frac{m^3-m}{hl^2}.\ee
Notice that the Hessian of $\eta$ is given by
\be\left[
\begin{array}{cc}
\eta_{ll} & \eta_{lm}\\
\eta_{lm} & \eta_{mm}
\end{array}\right]=-\frac{2}{l}\left[
\begin{array}{cc}
\frac{m^2+1}{l^2}  & -\frac{m}{l}\\
-\frac{m}{l} & 1
\end{array}\right].
\ee
Since $v-u>0$, we have $l<0$ and therefore the Hessian will be positive definite. Hence $\eta$ is a convex entropy. Next multiply system~\eqref{SV:lm viscous} by $(\eta_l,\eta_m)$ to get
\be \eta_t+q_x=C(x,t)+\mu(\eta_l\partial_{xx}l+\eta_m\partial_{xx}m)\ee
where $C(x,t)$ lies in a bounded set of $L^\infty(\Omega)$. Write $\eta_x=\eta_l\partial_xl+\eta_m\partial_xm$ and 
$$\eta_{xx}=\eta_{ll}(\partial_x l)^2+2\eta_{lm}\partial_x l\partial_x m+\eta_{mm}(\partial_x m)^2+\eta_l\partial_{xx} l+\eta_m\partial_{xx} m$$
so that
\be\label{SH: eta q}  \eta_t+q_x=C(x,t)-\mu(\eta_{ll}(\partial_xl)^2+2\eta_{lm}\partial_x l\partial_x m+\eta_{mm}(\partial_x m)^2   )+\mu\eta_{xx}\ee
Let $V$ be a compact subset of $\Omega$ and $\chi$ be a $C^\infty$ function with compact support in $\Omega$ and $\chi|_{V}=1$. Multiply~\eqref{SH: eta q} by $\chi$ and integrate over $\Omega$ to see
\begin{align}
\int_{\Omega} -\chi_t\eta-\chi_x q\,dtdx&=\int_{\Omega} \chi Cdtdx+\mu\int_{\Omega}\chi_{xx}\eta dtdx\nonumber\\
&-\mu\int_{\Omega}\chi(\eta_{ll}(\partial_xl)^2+2\eta_{lm}\partial_x l\partial_x m+\eta_{mm}(\partial_x m)^2   )dtdx.
\end{align}
Thus from the convexity of $\eta$ we have that 
\be \sqrt{\mu}\,\partial_x l,\qquad \sqrt{\mu}\,\partial_x m\ee
belong to a bounded subset of $L^2_{loc}(\Omega)$. Next compute
\begin{align} \mu\left|\int_{\Omega}\chi\partial_{xx}l \,dxdt\right|&=\mu\left| \int_{\Omega}\chi_x\partial_xl \,dx dt\right|\nonumber\\
&\le \sqrt{\mu}\left(\int_{supp(\chi)}\mu (\partial_x l)^2\,dxdt  \right)^{1/2}     \left(\int_{\Omega} (\chi_x)^2 dxdt   \right)^{1/2}  \nonumber\\
&\to 0\qquad\text{as}\,\mu\to0+.
\end{align}
Thus $\mu\partial_{xx} l\to0$ weakly in $L^2_{loc}(\Omega)$ and hence strongly in $H^{-1}_{loc}(\Omega)$. A similar statement holds for $\mu\partial_{xx} m$.

Hence we have proven

\begin{lemma}
 Under the assumptions of Lemma~\ref{Sbound: lemma}, for $u$, $v$ satisfying~\eqref{Sbound: lbd} at $t=T_1$, $T_1$ sufficiently large, we have that the sequences~\eqref{SH: seq} parametrized by $\mu$ lie in compact subsets of $H^{-1}_{loc}(\Omega)$.
\end{lemma}

\section{Compensated Compactness}\label{SCC}
In this section we use the a priori estimates of Section~\ref{SH} and the compensated compactness framework to pass to the limit as $\mu\to0+$ for our viscous system~\eqref{SV:lm viscous}. Notice the results of Section~\ref{SH} show that for initial data at $t=T_1$ , $T_1$ sufficiently large, which satisfy~\eqref{Sbound: lbd} we have
\be\label{A1} \left|(l^\mu,m^\mu,n^\mu  \right|\le A\qquad \text{for } (x,t)\in\Omega\ee
where $A$ is a constant independent of $\mu$ and
\be\label{A2}  \partial_t l^\mu-\partial_x(\frac{m^\mu}{h}),\quad \partial_t m^\mu-\partial_x(\frac{n^\mu}{h})\ee
confined in  a compact subset of $H^{-1}_{loc}(\Omega)$. Since~\eqref{A1}--\eqref{A2} are satisfied, Theorem $4.1$ of Chen, Slemrod and Wang~\cite{CSW} is applicable. We quote that result here.

\begin{theorem}[Compensated Compactness Framework~\cite{CSW}]
Let a sequence $(l^\mu, m^\mu, n^\mu)$ satisfy~\eqref{SV:lm viscous} and properties~\eqref{A1}--\eqref{A2}. There exists a subsequence, still labeled $(l^\mu, m^\mu, n^\mu)$ that converges weak* in $L^\infty(\Omega)$ to $(\tilde{l}, \tilde{m},\tilde{n})$ as $\mu\to0$ such that
\begin{enumerate}
\item[(i)] $|(\tilde{l}, \tilde{m},\tilde{n})(x,t)|\le A$\quad a.e. in $\Omega$
\item[(ii)] the Gauss equation~\eqref{SV:rel viscous} is weakly continuous with respect to the sequence $(l^\mu, m^\mu, n^\mu)$ that converges weak* in $L^\infty(\Omega)$ to $(\tilde{l}, \tilde{m},\tilde{n})$
\item[(iii)] the Codazzi equations hold for $(\tilde{l}, \tilde{m},\tilde{n})$.
\end{enumerate}
Specifically the limit $(\tilde{l}, \tilde{m},\tilde{n})$ is a bounded weak solution of the Gauss--Codazzi system in the domain $\Omega$.
\end{theorem}

The main result of this article is:

\begin{theorem}
Let $(\mathcal{M},g)$ be a geodesically complete simply connected smooth two
dimensional Riemannian manifold with non-positive Gauss curvature $K$ and a metric of the form~\eqref{S2 g}.
Assume that $h\doteq B$ and $k^*=|K|$ are independent of $x$ satisfying~\eqref{SV: h} and $k^*$ is given by~\eqref{SV: k}. There exists ${\bf y}\in W^{2,\infty}_{loc}(\Omega)$ satisfying the embedding equations of Theorem~\ref{Mardare}.
\end{theorem} 
\bp
Apply Mardare's theorem (see Theorem~\ref{Mardare}).
\ep

\section{Weaker decay and preservation of strict hyperbolicity}\label{S7}
Immediate inspection of the proof of Lemma~\ref{Sbound: lemma} shows that the bounds
\be\label{S7: uv} u\le -e^{-t}\psi_0,\qquad v\ge e^{-t}\psi_0\ee
for $t$ sufficiently large only follow from inequality~\eqref{Sbound: sign}. Hence a natural question is whether we can produce $k^*$, $h$ with weaker decay than given by~\eqref{SV: k} and still satisfy~\eqref{Sbound: sign} as well as $k^*$, $sk^*\in L^1[0,\infty)$. In fact the answer is yes as provided in the following example. Take
\be\label{S7: k} k^*=\frac{1}{(3+t)^2(\ln(3+t))^p},\qquad t>0,\,p> 1.\ee
Note that $k^*$, $t\,k^*(t)$ are in $L^1[0,\infty)$.

A direct computation shows
$$\frac{{k^*}'(t)}{k^*(t)}=-\frac{2}{3+t}-\frac{p}{(3+t)\ln(3+t)}.$$
Also recall from~\eqref{SV:la}--\eqref{SV:lb} that
$$ \frac{h'(t)}{h(t)}\ge\frac{\int_0^t k^*(s)ds}{1+C_1 t}$$
where $C_1$ is as given in Han-hong~\cite[Lemma 10.2.3]{HH} by the expression
$$ C_1=\int_0^\infty k^*(s)ds \exp\left\{\int_0^\infty sk^*(s)ds\right\}.$$
Hence~\eqref{Sbound: sign} will be satisfied if
\be\frac{\int_0^ tk^*(s) ds}{1+C_1 t}+\frac{1}{4}\left(-\frac{2}{3+t}-\frac{p}{(3+t)\ln(3+t)}\right)>0,\ee
or alternatively
\be\int_0^t k^*(s)ds-\frac{(1+C_1 t)}{4}\left(\frac{2}{3+t}+\frac{p}{(3+t)\ln(3+t)}\right)>0.\ee
Thus for large $t$ it suffices that
\be\label{S7: inq}\int_0^\infty k^*(s)ds>\frac{C_1}{2}.\ee
From the definition of $C_1$ we see~\eqref{S7: inq} will be satisfied when
\be 2>\exp\left\{\int_0^\infty k^*(s)ds\right\}.\ee
An easy estimate shows
\begin{align}
\int_0^\infty sk^*(s)ds&\le \int_0^\infty (s+3) k^*(s)ds\nonumber\\
&\le\int_0^\infty\frac{ds}{(s+3)\ln(3+s)^p}\nonumber\\
&=(\frac{1}{p-1})(\ln 3)^{1-p}.
\end{align}
Hence~\eqref{S7: inq} will be satisfied if
$$ 2>\exp\left\{\frac{1}{(p-1)(\ln 3)^{p-1}}\right\}.$$
As $\ln 3=1.0986\dots$ we want
\be\label{S7: inq2} 2>\exp\left\{\frac{1}{(p-1)(1.09)^{p-1}}\right\}\ee
and as $p\to\infty$ the right-hand side of~\eqref{S7: inq2} approaches $1$ and the inequality~\eqref{S7: inq2} is satisfied for $p$ large enough.

We summarize our observations in the following theorem.

\begin{theorem}
Assume that $h$ and $k^*$ are independent of $x$ satisfying~\eqref{SV: h} and $k^*$ is as given by~\eqref{S7: k}. Then if~\eqref{S7: uv} is satisfied at $t=T_1$, $T_1$ sufficiently large, then~\eqref{S7: uv} is satisfied for all $t\ge T_1$. Hence strict hyperbolicity of~\eqref{SV:lm}--\eqref{SV:rel} would not be lost and furthermore $\frac{-2 e^t}{\psi_0}<l<0$.

\end{theorem}

\subsection*{Acknowledgment}  
Christoforou was partially supported by the Start-Up fund 2011-2013 from University of Cyprus. Part of this work was carried out at Weizmann Institute in February of 2014 and Christoforou would like to thank the Department of Computer Science and Applied Mathematics at Weizmann Institute for the invitation and the hospitality. M.~Slemrod was supported in part by the
Simons Foundation Collaborative Research Grant $232531$.


\end{document}